\documentclass[11pt,a4paper]{article}
\usepackage{mathrsfs}
\usepackage{amsfonts}
\usepackage{amssymb}
\usepackage{amsmath}

\renewcommand{\thefootnote}{\fnsymbol{footnote}}

\begin{document}
\pagestyle{plain}
\renewcommand{\thefootnote}{\fnsymbol{footnote}}
\title{\Large \bf Landau's theorems for certain biharmonic mappings}
\author{\small \ Ming-Sheng Liu\thanks{Corresponding author.} \,  and Zhi-Wen Liu\\
{\small \ (School of Mathematical Sciences, South China Normal University}\\
{\small \ Guangzhou 510631, Guangdong, People's Republic of China)}\\
{\small \ E-mail:liumsh@scnu.edu.cn, liuziwen1985@163.com}\\{\small \  Yu-Can Zhu}\\
 { \small (Department of Mathematics, Fuzhou University,}\\
 { \small Fuzhou 350108, Fujian, P.R.China)}\\{ \small
 E-mail:zhuyucan@fzu.edu.cn}}
\date{}
\maketitle \footnote[0]{This research is partly supported by  the
Natural Science Foundation of  Fujian Province, China
(No.2009J01007) and the Education Commission Foundation of Fujian
Province, China (No.JA08013).}

\begin{center}

{\bf Published in Acta Mathematica Sinica, Chinese Series, 2011, Vol.54, No.1, 69-80.}

\end{center}

{\bf Abstract:}\ Let $f(z)=h(z)+\overline{g(z)}$ be a harmonic
mapping of the unit disk $U$. In this paper, the sharp coefficient
estimates for bounded planar harmonic mappings are established, the
sharp coefficient estimates for normalized planar harmonic mappings
with $|h(z)|+|g(z)|\leq M$ are also provided. As their applications,
Landau's theorems for certain biharmonic mappings are provided,
which improve and refine the related results of earlier authors.



\section{Introduction and preliminaries}
\mbox{}\indent Suppose that $f(z)=u(z)+iv(z)$ is a two times
continuously differentiable complex-valued function in a domain
$D\subseteq\mathbb{C}$. Then $f(z)$ is harmonic function in a domain
$D\subseteq\mathbb{C}$ if and only if $f(z)$ satisfies the following
harmonic equation
\begin{equation*}
\bigtriangleup f= 4f_{z\overline{z}}=\frac{\partial ^{2}f}{\partial
x^{2}}+\frac{\partial ^{2}f}{\partial y^{2}}=0,\; z=x+iy\in D,
\end{equation*}
where we use the common notations for its formal derivatives:
\begin{equation*}
f_{z}=\frac{1}{2}(f_{x}-if_{y}),\;
f_{\overline{z}}=\frac{1}{2}(f_{x}+if_{y}).
\end{equation*}

A four times continuously differentiable complex-valued function
$F(z)=U(z)+iV(z)$ is said to be a biharmonic in a domain
$D\subseteq\mathbb{C}$ if and only if $\bigtriangleup F$ is harmonic
in the domain $D$, i.e., $F(z)$ satisfies the following biharmonic
equation
$$
\bigtriangleup^2 F=\bigtriangleup (\bigtriangleup  F)=0, \;
z=x+iy\in D.
$$

Biharmonic functions arise in many physical situations, particularly
in fluid dynamics and elasticity problems, and have many important
applications in engineering(see \cite{amk} for the details).

Notice that the composition $f\circ \phi$ of a harmonic function $f$
with holomorphic function $\phi$ is harmonic, while this is not true
when $f$ is biharmonic. Without loss of generality, we consider the
class of biharmonic mappings defined in the open unit disk
$U=\{z\in\mathbb{C}:|z|<1\}$.

For such function $f$, we define
\begin{equation*}
\Lambda_{f}(z)=\max_{0\leq \theta \leq
2\pi}|f_{z}(z)+{e}^{-2i\theta}f_{\overline{z}}(z)|=|f_{z}(z)|+|f_{\overline{z}}(z)|,
\end{equation*}
and
\begin{equation*}
\lambda_{f}(z)=\min_{0\leq \theta \leq
2\pi}|f_{z}(z)+{e}^{-2i\theta}f_{\overline{z}}(z)|=||f_{z}(z)|-|f_{\overline{z}}(z)||.
\end{equation*}

We denote the Jacobian of $f$ by $J_{f}$, Lewy\cite{hle} showed that
a harmonic mapping $f(z)$ is locally univalent in a domain $D$ if
and only if $J_{f}(z)\neq 0$ for any $z\in D$. Of course, local
univalence of $f$ does not imply global univalence in a domain $D$.
If $D$ is simply connected, and $f(z)$ is a harmonic mapping in $D$,
then $f(z)$ can be written as $f=f_{1}+\overline{f_{2}}$ with
$f(0)=f_{1}(0)$, $f_{1}$ and $f_{2}$ are analytic on $D$. Thus, we
have
\begin{equation*}
J_{f}(z)=|f_{z}(z)|^{2}-|f_{\overline{z}}(z)|^{2}=|f'_{1}(z)|^{2}-|f'_{2}(z)|^{2}.
\end{equation*}

It is known \cite{amk} that a mapping $F$ is biharmonic in a simply
connected domain $D$ if and only if $F$ has the following
representation:
\begin{equation}
F(z)=|z|^{2}g(z)+h(z),\,\,\,z\in D,\label{l11}
\end{equation}
where $g(z)$ and $h(z)$ are complex-valued harmonic mappings in $D$.
It is well known that $g(z)$ and $h(z)$ can be expressed as
\begin{equation}
g(z)=g_{1}(z)+\overline{g_{2}(z)},\,\,\,z\in D,\label{l12}
\end{equation}
\begin{equation}
h(z)=h_{1}(z)+\overline{h_{2}(z)},\,\,\,z\in D,\label{l13}
\end{equation}
where $g_{1}(z)$, $g_{2}(z)$, $h_{1}(z)$ and $h_{2}(z)$ are analytic
in $D$ (see \cite{ys} for the details ).

The classical Landau theorem concerns determining the possibly
largest schlicht disk for the properly normalized bounded analytic
functions. It states that if $f$ is an analytic function on the unit
disk $U$ with $f(0)=f'(0)-1=0$ and $|f(z)|<M$ for $z\in U$, then $f$
is univalent in the disk $|z|<r_{0}$ with
$r_{0}=1/(M+\sqrt{M^{2}-1})$, and $f(|z|<r_{0})$ contains a disk
$|w|<R_{0}$ with $R_{0}=Mr_{0}^{2}$. This result is sharp, with the
extremal function $f(z)=Mz(1-Mz)/(M-z)$ (see \cite{chh}).

For harmonic mappings in $U$, under suitable restriction, Chen et
al. \cite{chh2} obtained two versions of Landau's
 theorems. These versions are different from each other by the
 normalization conditions assumed for the bounded harmonic mapping.
 However, their results are not sharp. Better estimates
 were given in \cite{mdm}, late in \cite{agr,lms1,lms2}.

For $r>0$, we let $U_{r}$ denote the disk with center at the origin
and radius $r$. Abdulhadi and Muhanna \cite{am2} obtained two
versions of Landau's theorems for biharmonic mappings; however,
their results are not sharp.

{\bf Theorem A}(Abdulhadi and Muhanna\cite{am2}) Let $f(z)=|z|^2
g(z)+h(z)$ be a biharmonic mapping of the unit disk $\mathbb{U}$, as
in (\ref{l11}), with $f(0)=h(0)=J_f(0)-1=0$ and $|g(z)|\leq M, \,
|h(z)|\leq M$ for $z\in\mathbb{U}$. Then there is a constant
$0<\rho_1<1$ so that $f$ is univalent in the disk
$\mathbb{U}_{\rho_1}$. In specific $\rho_1$ satisfies the following
equation
\begin{equation}
\frac{\pi}{4M} -2\rho_1 M -\frac{2M
\rho_1^2}{(1-\rho_1)^2}-2M\cdot\frac{2\rho_1-\rho_1^2}{(1-\rho_1)^2}=0,
\label{a1}
\end{equation}
and $f(\mathbb{U}_{\rho_1})$ contains a schlicht disk
$\mathbb{U}_{R_1}$ with
\begin{equation}
R_1=\frac{\pi}{4M}\rho_1- 2M\frac{\rho_1^3+\rho_1^2}{1-\rho_1} \, .
\label{a2}
\end{equation}

{\bf Theorem B}(Abdulhadi and Muhanna\cite{am2}) Let $g(z)$ be
harmonic in the unit disk $\mathbb{U}$, with $g(0)=J_g(0)-1=0$ and
$|g(z)|\leq M$ for $z\in\mathbb{U}$. Then $f(z)=|z|^2 g(z)$ is
univalent in the disk $\mathbb{U}_{\rho_2}$ with
\begin{equation}
\rho_2=\frac{\pi}{\pi+16M^2+2M\sqrt{2\pi+64M^2}}, \label{a01}
\end{equation}
and $f(\mathbb{U}_{\rho_2})$ contains a schlicht disk
$\mathbb{U}_{R_2}$ with
\begin{equation}
R_2=\frac{\pi}{4M}\rho_2^3- 2M\frac{\rho_2^4}{1-\rho_2} \, .
\label{a02}
\end{equation}

Recently, Liu \cite{lms} established the following better
coefficient estimates for bounded and normalized harmonic mappings.

{\bf Theorem C} (Liu \cite{lms}) Suppose that
$f(z)=h(z)+\overline{g(z)}$ is a harmonic mapping of the unit disk
$\mathbb{U}$ with $h(z)=\sum_{n=1}^{\infty}a_nz^n$ and
$g(z)=\sum_{n=1}^{\infty}b_nz^n$ for $z\in\mathbb{U}$.

(1)\ If $J_f(0)=1$ and $|f(z)|<M$, then
\begin{equation}
|a_n|+|b_n|\leq\sqrt{2M^2-2},\quad n=2, 3, \cdots ,\label{c02}
\end{equation}
and
\begin{equation}
\lambda_f(0)\geq\lambda_0(M)=\left\{
\begin{array}{lll}
\frac{\sqrt{2}}{\sqrt{M^2-1}+\sqrt{M^2+1}},&& \mbox{ if } 1\leq
M\leq M_0=
\frac{\pi}{2\sqrt[4]{2\pi^2-16}},\\\\
\frac{\pi}{4M}, && \mbox{ if } M>M_0=
\frac{\pi}{2\sqrt[4]{2\pi^2-16}}\approx1.1296 \,  .
\end{array}
\right. \label{l16}
\end{equation}

(2)\ If $\lambda_f(0)=1$ and $|f(z)|<M$, then the inequalities
(\ref{c02}) also hold.

By applying Theorem C, Liu improved Theorem A and Theorem B, and
obtained some completely new results. Wang  et al. \cite{cpw}
obtained two versions of Landau's
 theorems for biharmonic mappings of the form $L(F)$, where $F$ is a biharmonic mapping,
 and $L=z\frac{\partial}{\partial z}-\overline{z}\frac{\partial}{\partial\overline{z}}$.

{\bf Theorem D} \,(Liu \cite{lms}) Let $F(z)=|z|^2 g(z)+h(z)$ be a
biharmonic mapping of the unit disk $U$, as in (\ref{l11}), with
$F(0)=h(0)=\lambda_F(0)-1=0$ and $|g(z)|\leq M_1, \, |h(z)|\leq M_2$
for $z\in U$. Then, $F$ is univalent in the disk $U_{\rho_3}$, and
$F(U_{\rho_3})$ contains a schlicht disk $U_{R_3}$, where $\rho_3$
is the minimum positive root of the following equation
\begin{equation*}
1-2r M_1 -2M_1\cdot\frac{r^2}{(1-r)^2}
-\sqrt{2M_2^2-2}\cdot\frac{2r-r^2}{(1-r)^2}=0,
\end{equation*}
and
\begin{equation*}
R_3=\rho_3-\frac{2M_1 \rho_3^3}{1-\rho_3}-\sqrt{2M_2^2-2}\cdot\frac{
\rho_3^2}{1-\rho_3} \, .
\end{equation*}

{\bf Theorem E} \,(Liu \cite{lms}) Let $F(z)=|z|^{2}g(z)+h(z)$ be a
biharmonic  mapping of the unit disk $U$, as in (\ref{l11}), with
$F(0)=h(0)=J_{F}(0)-1=0$, and $|g(z)|\leq M_{1}$, $|h(z)|\leq M_{2}$
for $z\in U$. Then $F$ is univalent in the disk $U_{\rho_{4}}$ , and
$F(U_{\rho_{4}})$ contains a schlicht disk $U_{R_{4}}$, where
$\rho_{4}$ is the minimum positive root of the following equation:
\begin{equation}
\lambda _{0}(M_{2})-2r
M_{1}-\frac{2M_{1}r^{2}}{(1-r)^{2}}-\sqrt{2M_{2}^{2}-2}\cdot
\frac{2r-r^{2}}{(1-r)^{2}}=0,\label{l14}
\end{equation}
and
\begin{equation}
R_{4}=\lambda _{0}(M_{2})\rho_{4}
-\frac{2M_{1}\rho_{4}^{3}}{1-\rho_{4}}-\sqrt{2M_{2}^{2}-2}\cdot
\frac{\rho_{4}^{2}}{1-\rho_{4}},\label{l15}
\end{equation}
where $\lambda _{0}(M)$ is defined by (\ref{l16}).

{\bf Theorem F} \,(Liu \cite{lms}) Let $g(z)$ be harmonic in the
unit disk $U$, with $g(0)=\lambda_{g}(0)-1=0$ and $|g(z)|\leq M$ for
$z\in U$. Then $F(z)=|z|^{2}g(z)$ is univalent in the disk
$U_{\rho_{5}}$, and $F(U_{\rho_{5}})$ contains a schlicht disk
$U_{R_{5}}$, where
\begin{eqnarray}
\rho_{5}=\frac{1}{1+2\sqrt{2M^{2}-2}+\sqrt{\sqrt{2M^{2}-2}+8(M^{2}-1)}},\label{l17}
\end{eqnarray}
and
\begin{eqnarray}
R_{5}=\left\{ \begin{array}{lll}
\rho_{5}^{3}-\sqrt{2M^{2}-2}\cdot\frac{\rho_{5}^{4}}{1-\rho_{5}},
\qquad &M>1,\\
1,\qquad &M=1.
\end{array} \right.\label{l18}
\end{eqnarray}
Above result is sharp when $M=1$.

{\bf Theorem G} \,(Liu \cite{lms})  Let $g(z)$ be harmonic in the
unit disk $\mathbb{U}$, with $g(0)=J_g(0)-1=0$ and $|g(z)|\leq M$
for $z\in\mathbb{U}$. Then $F(z)=|z|^2 g(z)$ is univalent in the
disk $\mathbb{U}_{\rho_6}$, and $F(\mathbb{U}_{\rho_6})$ contains a
schlicht disk $\mathbb{U}_{R_6}$, where
\begin{equation}
\rho_6=\frac{\lambda_0(M)}{\lambda_0(M)+2\sqrt{2M^2-2}+\sqrt{\lambda_0(M)
\sqrt{2M^2-2}+8(M^2-1)}}, \label{b15}
\end{equation}
and
\begin{equation}
R_6=\left\{
\begin{array}{lll}
\lambda_0(M) \, \rho_6^3-\sqrt{2M^2-2}\cdot\frac{\rho_6^4}{1-\rho_6}
\, , &\mbox{ if }M>1,\\
1,&\mbox{ if }M=1,
\end{array}
\right.\label{b16}
\end{equation}
where $\lambda_0(M)$ is defined by (\ref{l16}). Above result is
sharp when $M=1$.

The coefficient estimates (\ref{c02}) for bounded and normalized
harmonic mappings in Theorem C are not sharp for $M>1$. Furthermore,
Liu\cite{lms} proposed the following conjecture:

{\bf Conjecture H} \,(Liu \cite{lms}) Suppose that
$f(z)=h(z)+\overline{g(z)}$ is a harmonic mapping of the unit disk
$U$ with $h(z)=\sum_{n=1}^\infty a_{n}z^{n}$ and
$g(z)=\sum_{n=1}^\infty b_{n}z^{n}$ are analytic on $U$, moreover
$|f(z)|\leq M$ for $z\in U$. If  $J_{f}(0)=1$ or $\lambda_{f}(0)=1$,
then
\begin{eqnarray*}
|a_{n}|+|b_{n}|\leq M-\frac{1}{M},\,\,\, n=2,3,\cdots,
\end{eqnarray*}
with the extremal functions $f_{n}(z)=Mz((1-Mz^{n-1})/(M-z^{n-1}))$
for $n=2,3,\cdots$.

In this paper, we first establish the sharp coefficient estimates
for bounded harmonic mappings (Lemma 2.1). Next, we prove that the
above conjecture is true for the subclass of harmonic mappings
$f(z)=h(z)+\overline{g(z)}$ which satisfy the following bounded
condition
\begin{eqnarray}
|h(z)|+|g(z)|\leq M.\label{l19}
\end{eqnarray}
Under this condition, we establish the sharp coefficient estimates
for normalized harmonic mappings(Lemma 2.3, Corollary 2.4). Then,
using these estimates, we verify several versions of Landau's
theorem (Theorems 2.6, 2.6$'$, 2.8, 2.8$'$, 2.10, and 2.10$'$,
Corollary 2.12, 2.12$'$), which refine and improve Theorem A$-$G. In
order to establish our main results, we recall the following lemma.

{\bf Lemma 1.1} \, (see, for example,p.35 in \cite{ig}) If
$F(z)=a_{0}+a_{1}z+\cdots+a_{n}z^{n}+\cdots$ is analytic and
$|F(z)|\leq1$ on $U$. Then $|a_{n}|\leq 1-|a_{0}|^{2}$ for
$n=1,2,\cdots$.

\setcounter{equation}{0}
\section{Main results}
\mbox{}\indent We first establish the following sharp coefficient
estimates for bounded harmonic mappings.

{\bf Lemma 2.1} \, Suppose that $f(z)=h(z)+\overline{g(z)}$ is a
harmonic mapping of the unit disk $U$ with
$h(z)=\sum\limits_{n=0}^{\infty} a_{n}z^{n}$ and
$g(z)=\sum\limits_{n=1}^{\infty }b_{n}z^{n}$. If $|f(z)|\leq M$ for
$z\in U$, then we have
\begin{eqnarray}
|a_n|+|b_n|\leq\frac{4}{\pi}M,  \, n=1, 2, \cdots.\label{a3}
\end{eqnarray}
The inequalities (\ref{a3}) are sharp for all $n=1,2,\cdots$.

{\bf Proof} \,  Fix $r\in (0, 1)$, we have
$$
f(re^{i\theta})=\sum\limits_{n=0}^{\infty}a_nr^ne^{in\theta}+\sum\limits_{n=1}^{\infty}\overline{b}_nr^ne^{-in\theta},
\, \theta\in [0, 2\pi).
$$
Hence we obtain
\begin{eqnarray}
a_nr^n=\frac{1}{2\pi}\int_0^{2\pi}f(re^{i\theta})e^{-in\theta}d\theta,\label{a5}
\,  n=0, 1, 2, \cdots,
\end{eqnarray}
and
\begin{eqnarray}
\overline{b}_nr^n=\frac{1}{2\pi}\int_0^{2\pi}f(re^{i\theta})e^{in\theta}d\theta,
\,  n=1, 2, \cdots.\label{a6}
\end{eqnarray}

For every $n$, we let $a_n=|a_n|e^{i\alpha_n},
b_n=|b_n|e^{i\beta_n}, \theta_n=\frac{\beta_n+\alpha_n}{2n}$. Since
$|\cos n\theta|$ is a periodic function with period $\pi$, from
(\ref{a5}), (\ref{a6}), we have
\begin{eqnarray*}
(|a_n|+|b_n|)r^n&=&\bigg|\frac{1}{2\pi}\int_0^{2\pi}f(re^{i\theta})
[e^{-i\alpha_n}e^{-in\theta}+e^{i\beta_n}e^{in\theta}]d\theta\bigg|\\
&\leq&\frac{1}{2\pi}\int_0^{2\pi}|f(re^{i\theta})||e^{-i\alpha_n}e^{-in\theta}+e^{i\beta_n}e^{in\theta}|d\theta\\
&=&\frac{1}{2\pi}\int_0^{2\pi}|f(re^{i\theta})||e^{-in\theta}+e^{i(\beta_n+\alpha_n)}e^{in\theta}|d\theta\\
&=&\frac{1}{2\pi}\int_0^{2\pi}|f(re^{i\theta})||e^{-in(\theta+\theta_n)}+e^{in(\theta+\theta_n)}|d\theta\\
&\leq&\frac{1}{2\pi}\int_0^{2\pi}2M|\cos
n(\theta+\theta_n)|d\theta\\
&=&\frac{M}{\pi}\int_{0}^{2\pi}|\cos n\theta|d\theta =\frac{M}{n\pi}\int_{0}^{2n\pi}|\cos x|dx \\
&=&\frac{4M}{\pi}.
\end{eqnarray*}

Setting $r\rightarrow 1^-$, we obtain
$|a_n|+|b_n|\leq\frac{4M}{\pi}$ for all $n=1,2,\cdots$.

In order to verify the sharpness of (\ref{a3}), we choose the
harmonic mapping
$$
F(z)=\mbox{Im}\bigg\{\frac{2M}{\pi}\log\frac{1+z}{1-z}\bigg\}=\frac{2M}{\pi}\arctan\frac{2y}{1-x^2-y^2},
\, z=x+iy\in U.
$$
Then we have $|F(z)|\leq M$ for $z\in U$.  By a direct calculation,
we obtain
$$
F(z)=\sum\limits_{n=0}^{\infty}A_nz^n+\overline{\sum\limits_{n=1}^{\infty}B_nz^n}=-\frac{2Mi}{\pi}z+\cdots+\overline{\frac{2Mi}{\pi}z+\cdots}.
$$
This implies $|A_1|+|B_1|=\frac{4M}{\pi}$. Let $m$ be an integer
number with $m>1$. Then the mapping
$$
F_m(z)=F(z^m)=\sum\limits_{n=m}^{\infty}C_nz^n+\overline{\sum\limits_{n=m}^{\infty}D_nz^n}
=-\frac{2Mi}{\pi}z^m+\cdots+\overline{\frac{2Mi}{\pi}z^m+\cdots}
$$
is  harmonic in $U$ with $|F_m(z)|\leq M$ for $z\in U$ such that
$|C_m|+|D_m|=\frac{4M}{\pi}$. This completes the proof.\hfill
$\square$

{\bf Remark 2.2} \, If
$M>M_0'=\frac{\pi}{\sqrt{\pi^2-8}}\approx2.2976$, then
$\sqrt{2M^2-2}>\frac{4M}{\pi}$.  From Lemma 2.1, we improve Theorem
C or Lemma 2.1 in \cite{lms2} for
$M>M_0'=\frac{\pi}{\sqrt{\pi^2-8}}\approx2.2976$.

Next, we establish the following sharp coefficient estimates for
normalized harmonic mappings which  satisfy the bounded condition
(\ref{l19}).

{\bf Lemma 2.3} \, Suppose that $f(z)=h(z)+\overline{g(z)}$ is a
harmonic mapping of the unit disk $U$ with $h(z)=\sum_{n=1}^\infty
a_{n}z^{n}$ and $g(z)=\sum_{n=1}^\infty b_{n}z^{n}$ are analytic on
$U$. If $|h(z)|+|g(z)|\leq M$ for $z\in U$, then we have $0\leq
\lambda_{f}(0)\leq M$, and
\begin{eqnarray}
|a_{n}|+|b_{n}|\leq
M-\frac{\lambda_{f}^{2}(0)}{M},\,\,\,n=2,3,\cdots.\label{l21}
\end{eqnarray}
Above estimates are sharp for all $n=2,3,\cdots$, with the extremal
functions $f_{a,\ n}(z)$ and $\overline{f_{a,\ n}(z)}$
\begin{eqnarray}
f_{a,\ n}(z)=Mz\cdot\frac{a-Mz^{n-1}}{M-az^{n-1}},\label{l22}
\end{eqnarray}
where $0\leq a\leq M$.

{\bf Proof} \, Fix $n\in \mathbb{N}-\{1\}=\{2,3,\ldots\}$, we choose
a real number $\alpha$ such that
$|a_{n}+e^{i\alpha}b_{n}|=|a_{n}|+|b_{n}|$, and set
\begin{eqnarray*}
l(z)=\frac{1}{M}[h(z)+e^{i\alpha}g(z)]=\sum_{k=1}^\infty
\frac{a_{k}+e^{i\alpha}b_{k}}{M}z^{k}.
\end{eqnarray*}

Since $h(z)$ and $g(z)$ are analytic and $|f(z)|\leq
|h(z)|+|g(z)|\leq M$ on $U$, we get that $l(z)$ is analytic and
$|l(z)|\leq (|h(z)|+|g(z)|)/M\leq 1$ on $U$. Notice $l(0)=0$, by
Schwarz lemma, we obtain that $|l(z)|\leq |z|$. Setting
\begin{eqnarray*}
F(z)=\frac{a_{1}+e^{i\alpha}b_{1}}{M}+\sum_{k=2}^\infty
\frac{a_{k}+e^{i\alpha}b_{k}}{M}z^{k-1},
\end{eqnarray*}
it follows that $F(z)$ is analytic and $|F(z)|\leq1$ on $U$. Since
\begin{eqnarray}
||a_{1}|-|b_{1}||=\lambda_{f}(0),\label{l23}
\end{eqnarray}
by Lemma 1.1 and (\ref{l23}), we obtain
\begin{eqnarray*}
\left|\frac{a_{k}+e^{i\alpha}b_{k}}{M}\right|\leq1-\left|\frac{a_{1}+e^{i\alpha}b_{1}}{M}\right|^{2}
\leq1-\frac{||a_{1}|-|b_{1}||^{2}}{M^{2}}=1-\frac{\lambda_{f}^{2}(0)}{M^{2}},\,\,\,
k=2,3,\cdots.
\end{eqnarray*}
Thus we have $0\leq \lambda_{f}(0)\leq M$. In particular, we have
\begin{eqnarray*}
|a_{n}|+|b_{n}|=|a_{n}+e^{i\alpha}b_{n}|\leq
M\bigg(1-\frac{\lambda_{f}^{2}(0)}{M^{2}}\bigg)=M-\frac{\lambda_{f}^{2}(0)}{M}.
\end{eqnarray*}

Finally, it is obvious that the equalities hold for all
$n=2,3,\cdots$ for the functions
$$f_{a,\
n}(z)=Mz\frac{a-Mz^{n-1}}{M-az^{n-1}}=az-(M-\frac{a^2}{M})z^n+\cdots$$
and
$$
\overline{f_{a,\
n}(z)}=a\overline{z}-(M-\frac{a^2}{M})\overline{z}^n+\cdots,
$$
respectively, where $a=\lambda_{f}(0)$. This completes the
proof.\hfill $\square$

Setting $\lambda_{f}(0)=1$ in Lemma 2.3, we get the following
corollary.

{\bf Corollary 2.4}\, Suppose that $f(z)=h(z)+\overline{g(z)}$ is a
harmonic mapping of the unit disk $U$ with $h(z)=\sum_{n=1}^\infty
a_{n}z^{n}$ and $g(z)=\sum_{n=1}^\infty b_{n}z^{n}$ are analytic on
$U$. If $|h(z)|+|g(z)|\leq M$ for $z\in U$, and $\lambda_{f}(0)=1$,
then $M\geq 1$, and
\begin{eqnarray*}
|a_{n}|+|b_{n}|\leq M-\frac{1}{M},\,\,\,n=2,3,\cdots.%\label{l22}
\end{eqnarray*}
Above estimates are sharp for all $n=2,3,\cdots$, with the extremal
functions $f_{1, n}(z)$ and $\overline{f_{1, n}(z)}$ defined by
(\ref{l22}).

{\bf Remark 2.5} \, Corollary 2.4 tell us that Conjecture H is true
for the subclass of harmonic mappings $f$ which satisfy the bounded
condition (\ref{l19}) and $\lambda_{f}(0)=1$. Setting $M=1$ in
Corollary 2.4, we get that $a_{n}=b_{n}=0$ for $n=2,3,\ldots$, thus
$f(z)=\alpha z$ or $f(z)=\alpha\overline{z}$, where $|\alpha|=1$.

Now, with the aid of Lemmas 2.1 and 2.3, we can improve Theorem D as
follows.

{\bf Theorem 2.6} \, Let $F(z)=|z|^{2}g(z)+h(z)$ be a biharmonic
mapping of the unit disk $U$, as in (\ref{l11})-(\ref{l13}), with
$F(0)=h(0)=\lambda_{F}(0)-1=0$, and $|g(z)|\leq M_{1}$, $|h(z)|\leq
 M_{2}$ for $z\in U$.
Then $F$ is univalent in the disk $U_{r_{1}}$, and
$F(U_{\sigma_{1}})$ contains a schlicht disk $U_{\sigma_{1}}$, where
$r_{1}$ is the minimum positive root of the following equation:
\begin{equation}
1-2r M_{1}-\frac{4 M_{1}r^{2}}{\pi (1-r)^{2}} -K(M_{2})\cdot
\frac{2r-r^{2}}{(1-r)^{2}}=0,\label{l25}
\end{equation}
and
\begin{equation}
\sigma_{1}=r_{1}-\frac{4M_{1}r_{1}^{3}}{\pi (1-r_{1})}-K(M_{2})\cdot
\frac{r_{1}^{2}}{1-r_{1}},\label{l26}
\end{equation}
where $K(M_{2})=\min\{\sqrt{2M_{2}^2-2}, \frac{4}{\pi}M_{2}\}$.

{\bf Proof} \, From (\ref{l12}) and (\ref{l13}), we have
\begin{equation*}
g(z)=g_{1}(z)+\overline{g_{2}(z)},\quad
h(z)=h_{1}(z)+\overline{h_{2}(z)},
\end{equation*}
where $g_{1}(z)=\sum_{n=0}^\infty a_{n}z^{n}$,
$g_{2}(z)=\sum_{n=1}^\infty b_{n}z^{n}$, $h_{1}(z)=\sum_{n=1}^\infty
c_{n}z^{n}$ and $h_{2}(z)=\sum_{n=1}^\infty d_{n}z^{n}$ are analytic
in $U$, then, we have
\begin{equation}
\lambda_{F}(0)=||c_{1}|-|d_{1}||=\lambda_{h}(0)=1.\label{l27}
\end{equation}

To prove the univalence of $F(z)$ in $U_{r_{1}}$, we adopt the
method used in \cite{am2,lms}. By means of the hypothesis of Theorem
2.5, by Theorem C and Lemmas 2.1, we have
\begin{eqnarray}
|a_{n}|+|b_{n}|\leq
\frac{4}{\pi}M_{1}(n=1,2,\cdots),\,\,\,|c_{n}|+|d_{n}|\leq K(M_{2})\
(n=2,3,\cdots),\label{l28}
\end{eqnarray}
where $K(M_{2})=\min\{\sqrt{2M_{2}^2-2}, \frac{4}{\pi}M_{2}\}$.

Thus, for $z_{1}\neq z_{2}$ in $U_{r}(0<r<r_{1})$, by (\ref{l27})
and (\ref{l28}), we have
\begin{eqnarray*}
|F(z_{1})-F(z_{2})|&=&\left|\int_{[z_{1},z_{2}]}F_{z}(z)dz+F_{\overline{z}}(z)d\overline{z}\right|\\
&\geq&\left|\int_{[z_{1},z_{2}]}h_{z}(0)dz+h_{\overline{z}}(0)d\overline{z}\right|-
\left|\int_{[z_{1},z_{2}]}g(z)(\overline{z}dz+zd\overline{z})\right|\\
&-&\left|\int_{[z_{1},z_{2}]}|z|^{2}(g_{1}'(z)dz+\overline{g_{2}'(z)}d\overline{z})\right|\\
&-&\left|\int_{[z_{1},z_{2}]}(h_{1}'(z)-h_{1}'(0))dz+(\overline{h_{2}'(z)}-\overline{h_{2}'(0)})d\overline{z}\right|\\
&\geq&|z_{1}-z_{2}|\left(\lambda_{h}(0)-2rM_{1}-\sum_{n=1}^\infty
(|a_{n}|+|b_{n}|)nr^{n+1}-\sum_{n=2}^\infty
(|c_{n}|+|d_{n}|)nr^{n-1}\right)\\
&\geq&|z_{1}-z_{2}|\left(1-2rM_{1}-\sum_{n=1}^\infty
\frac{4}{\pi}M_{1}nr^{n+1}-\sum_{n=2}^\infty
K(M_{2})nr^{n-1}\right)\\
&\geq&|z_{1}-z_{2}|\left(1-2rM_{1}-\frac{4M_{1}r^{2}}{\pi(1-r)^{2}}-
K(M_{2})\cdot\frac{2r-r^{2}}{(1-r)^{2}}\right)>0,
\end{eqnarray*}
this implies $F(z_{1})\neq F(z_{2})$.

Notice that $F(0)=0$, for any $z^{'}=r_{1}e^{i\theta}\in \partial
U_{r_{1}}$, by (\ref{l28}), we have
\begin{eqnarray*}
|F(z^{'})|&\geq&
|c_{1}z^{'}+\overline{d}_{1}\overline{z^{'}}|-r_{1}^{2}\left|\sum_{n=1}^\infty
a_{n}z^{'n}+\overline{b}_{n}\overline{z^{'}}^{n}\right|-\left|\sum_{n=2}^\infty
c_{n}z^{'n}+\overline{d}_{n}\overline{z^{'}}^{n}\right|\\
&\geq&\lambda_{h}(0)r_{1}-r_{1}^{2}\sum_{n=1}^\infty
(|a_{n}|+|b_{n}|)r_{1}^{n}-\sum_{n=2}^\infty
(|c_{n}|+|d_{n}|)r_{1}^{n}\\
&\geq&r_{1}-r_{1}^{2}\sum_{n=1}^\infty\frac{4}{\pi}M_{1}r_{1}^{n}-\sum_{n=2}^\infty K(M_{2})r_{1}^{n}\\
&=&r_{1}-\frac{4}{\pi}M_{1}\frac{r_{1}^{3}}{1-r_{1}}-K(M_{2})\frac{r_{1}^{2}}{1-r_{1}}
=\sigma_1.
\end{eqnarray*}

Hence, $F(z)$ is univalent on $U_{r_{1}}$ and $F(U_{r_{1}})$
contains the disk $U_{\sigma_{1}}$, where $r_{1}$ is defined by
(\ref{l25}) and $\sigma_{1}$ is defined by (\ref{l26}). This
completes the proof of Theorem 2.6.\hfill $\square$

With the aid of Lemma 2.1 and Corollary 2.4, applying the similar
method as in our proof of Theorem 2.6, we get a version of Landau's
theorem as follows.

{\bf Theorem 2.6$'$} \, Let $F(z)=|z|^{2}g(z)+h(z)$ be a biharmonic
mapping of the unit disk $U$, as in (\ref{l11})-(\ref{l13}), with
$F(0)=h(0)=\lambda_{F}(0)-1=0$, and $|g(z)|\leq M_{1}$,
$|h_{1}(z)|+|h_{2}(z)|\leq M_{2}$ for $z\in U$. Then $F$ is
univalent in the disk $U_{r_{1}'}$, and $F(U_{\sigma_{1}'})$
contains a schlicht disk $U_{\sigma_{1}'}$, where $r_{1}'$ is the
minimum positive root of the following equation:
\begin{equation*}
1-2r M_{1}-\frac{4 M_{1}r^{2}}{\pi (1-r)^{2}}
-(M_{2}-\frac{1}{M_{2}})\cdot \frac{2r-r^{2}}{(1-r)^{2}}=0,
\end{equation*}
and
\begin{equation*}
\sigma_{1}=r_{1}'-\frac{4M_{1}r_{1}'^{3}}{\pi
(1-r_{1}')}-(M_{2}-\frac{1}{M_{2}})\cdot
\frac{r_{1}'^{2}}{1-r_{1}'}.
\end{equation*}

{\bf Remark 2.7} \, Notice that $\frac{4}{\pi}M_1<2M_1$,
$M_2-\frac{1}{M_{2}}<\frac{4}{\pi}M_2\leq\sqrt{2M_2^2-2}$ for
$M_2\geq M_0'=\frac{\pi}{\sqrt{\pi^2-8}}\approx2.2976$, it is easy
to verify that
$$
r_1'>r_1>\rho_3,\quad\sigma_{1}'>\sigma_{1}>R_3.
$$

With the aid of  Theorem C and Lemma 2.1, applying the similar
method as in our proof of Theorem 2.6, we can improve Theorem E as
follows.

{\bf Theorem 2.8} \, Let $F(z)=|z|^{2}g(z)+h(z)$ be a biharmonic
mapping of the unit disk $U$, as in (\ref{l11}), with
$F(0)=h(0)=J_{F}(0)-1=0$, and $|g(z)|\leq M_{1}$, $|h(z)|\leq M_{2}$
for $z\in U$, as in (\ref{l12}) and (\ref{l13}). Then $F$ is
univalent in the disk $U_{r_{2}}$ , and $F(U_{\sigma_{2}})$ contains
a schlicht disk $U_{\sigma_{2}}$, where $r_{2}$ is the minimum
positive root of the following equation:
\begin{equation}
\lambda_{0}(M_{2})-2rM_{1}-\frac{4M_{1}r^{2}}{\pi (1-r)^{2}}-
K(M_{2})\cdot\frac{2r-r^{2}}{(1-r)^{2}}=0,%\label{l29}
\end{equation}
and
\begin{equation}
\sigma_{2}=\lambda_{0}(M_{2})r_{2}-\frac{4M_1\ r_{2}^{3}}{\pi (1-r_{2})}-K(M_{2})\cdot\frac{r_{2}^{2}}{1-r_{2}},%\label{l210}
\end{equation}
where $\lambda _{0}(M)$ is defined by (\ref{l16}).

{\bf Proof} \, By the proof of Theorem 2.6 and the hypothesis of
Theorem 2.8, we have
\begin{equation*}
J_{F}(0)=|c_{1}|^{2}-|d_{1}|^{2}=J_{h}(0)=1.
\end{equation*}

By (\ref{l16}), we get that
\begin{equation}
\lambda_{h}(0)\geq\lambda_{0}(M_{2}).\label{l211}
\end{equation}

By means of Lemmas 2.1 and Theorem C, applying the similar method as
in our proof of Theorem 2.6, we can completes the proof of Theorem
2.8. \hfill $\square$

With the aid of Lemma 2.3 and (\ref{l16}), we can get the second
version of Landau's theorem.

{\bf Theorem 2.8$'$} \, Let $F(z)=|z|^{2}g(z)+h(z)$ be a biharmonic
mapping of the unit disk $U$, as in (\ref{l11})-(\ref{l13}), with
$F(0)=h(0)=J_{F}(0)-1=0$, and $|g(z)|\leq M_{1}$,
$|h_{1}(z)|+|h_{2}(z)|\leq M_{2}$ for $z\in U$. Then $F$ is
univalent in the disk $U_{r_{2}'}$ , and $F(U_{\sigma_{2}'})$
contains a schlicht disk $U_{\sigma_{2}'}$, where $r_{2}'$ is the
minimum positive root of the following equation:
\begin{equation}
\lambda_{0}(M_{2})-2rM_{1}-\frac{4M_{1}r^{2}}{\pi (1-r)^{2}}-
(M_{2}-\frac{\lambda_{0}^{2}(M_{2})}{M_{2}})\cdot\frac{2r-r^{2}}{(1-r)^{2}}=0,\label{l29}
\end{equation}
and
\begin{equation}
\sigma_{2}'=\lambda_{0}(M_{2})r_{2}'-\frac{4M_1\ r_{2}'^{3}}{\pi
(1-r_{2}')}-(M_{2}-
\frac{\lambda_{0}^{2}(M_{2})}{M_{2}})\cdot\frac{r_{2}'^{2}}{1-r_{2}'},\label{l210}
\end{equation}
where $\lambda _{0}(M)$ is defined by (\ref{l16}).

{\bf Remark 2.9} \, Notice that $\frac{4}{\pi}M_1<2M_1$,
$M_2-\frac{\lambda_{0}^{2}(M_{2})}{M_{2}}<\frac{4}{\pi}M_2\leq\sqrt{2M_2^2-2}$
for $M_2\geq M_0'=\frac{\pi}{\sqrt{\pi^2-8}}\approx2.2976$, it is
easy to verify that
$$
r_2'>r_2>\rho_4,\quad\sigma_{2}'>\sigma_{2}>R_4,
$$
for $M_2>1$.

The next theorem is different. Because, when $h=0$, the Jacobian
$J_{F}(0)=0$ and hence we assume that $\lambda_{g}(0)=1$ instead.
With the aid of  Theorem C and Lemma 2.1, we can improve Theorem F
as follows.

{\bf Theorem 2.10} \, Let $g(z)$ be harmonic in the unit disk $U$,
with $g(0)=\lambda_{g}(0)-1=0$ and $|g(z)|\leq  M$ for $z\in U$, as
in (\ref{l12}). Then
 $F(z)=|z|^{2}g(z)$ is univalent in the disk $U_{r_{3}}$, and
$F(U_{r_{3}})$ contains a schlicht disk $U_{\sigma_{3}}$, with
\begin{equation}
K(M)=\min\{\sqrt{2M^2-2},\ \frac{4M}{\pi}\},\label{l220}
\end{equation}
and
\begin{equation*}
r_{3}=\frac{1}{1+2K(M)+\sqrt{K(M)+4K(M)^2}},\label{l213}
\end{equation*}
and
\begin{eqnarray*}
\sigma_{3}= \left\{ \begin{array}{lll}
r_{3}^{3}-K(M)\cdot\frac{r_{3}^{4}}{1-r_{3}}, &M>1,\\
1 ,\qquad      & M=1.
\end{array} \right.\label{l214}
\end{eqnarray*}
Above result is sharp for $M=1$.

{\bf Proof.} If $F(z)=|z|^{2}g(z)$ satisfies the hypothesis of
Theorem 2.10, where
\begin{equation*}
g(z)=g_{1}(z)+\overline{g_{2}(z)}=\sum_{n=1}^\infty
a_{n}z^{n}+\overline{\sum_{n=1}^\infty b_{n}z^{n}}
\end{equation*}
is harmonic in $U$, then
\begin{equation}
||a_{1}|-|b_{1}||=\lambda_{g}(0)=1.\label{l216}
\end{equation}

By Theorem C and Lemma 2.1, we have
\begin{equation}
|a_{n}|+|b_{n}|\leq K(M) \, (n=2,3,\cdots).\label{l215}
\end{equation}

For $z_{1}\neq z_{2}$ in $U_{r}(0<r<r_{3})$, we have
\begin{eqnarray*}
F(z_{1})-F(z_{2})=\int_{[z_{1},z_{2}]}F_{z}(z)dz+F_{\overline{z}}(z)d\overline{z}
=\int_{[z_{1},z_{2}]}(\overline{z}g+|z|^{2}g'_{1})dz+(zg+|z|^{2}\overline{g'_{2}})d\overline{z}
\end{eqnarray*}
where $[z_{1},z_{2}]$ is the line segment from $z_{1}$ to $z_{2}$,
$z=(1-t)z_{1}+tz_{2}$ and $t\in[0,1]$. Notice that $r_{3}$ is the
minimum positive root of the following equation:
$$
1- K(M) \cdot\frac{4r-3r^{2}}{(1-r)^{2}}=0.
$$
Applying the same method used in \cite{am2,lms}, by (\ref{l215}) and
(\ref{l216}), we have
\begin{eqnarray*}
\left|\frac{F(z_{1})-F(z_{2})}{z_{1}-z_{2}}\right|&=&\frac{1}{|z_{1}-z_{2}|}
\left|\int_{[z_{1},z_{2}]}(\overline{z}g+|z|^{2}g'_{1})dz+
(zg+|z|^{2}\overline{g'_{2}})d\overline{z}\right|\\
&\geq&\int_{[z_{1},z_{2}]}|z|^{2}dt\cdot\left(||a_{1}|-|b_{1}||-2\sum_{n=2}^\infty
(|a_{n}|+|b_{n}|)r^{n-1}-\sum_{n=2}^\infty
(|a_{n}|+|b_{n}|)nr^{n-1}\right)\\
&\geq&\int_{[z_{1},z_{2}]}|z|^{2}dt\cdot\left(1-2\sum_{n=2}^\infty
 K(M) r^{n-1}-\sum_{n=2}^\infty
 K(M) nr^{n-1}\right)\\
&=&\int_{[z_{1},z_{2}]}|z|^{2}dt\cdot\left(1-2
 K(M) \cdot\frac{r}{1-r}-
 K(M) \cdot\frac{2r-r^{2}}{(1-r)^{2}}\right)\\
&=&\int_{[z_{1},z_{2}]}|z|^{2}dt\cdot\left(1-
 K(M) \cdot\frac{4r-3r^{2}}{(1-r)^{2}}\right)>0\\
\end{eqnarray*}
this implies $F(z_{1})\neq F(z_{2})$ for $z_{1}\neq z_{2}$ in
$U_{r}(0<r<r_{3})$, where $|z|=|(1-t)z_{1}+tz_{2}|$. Hence
$F(z)=|z|^{2}g(z)$ is univalent in the disk $U_{r_{3}}$.

From Theorem C, we get $M\geq1$. When $M>1$, by means of
(\ref{l216}), for $|z|=r_{3}$, we have
\begin{eqnarray*}
|F(z)|&=&r^{2}_{3}\left|\sum_{n=1}^\infty
a_{n}z^{n}+\overline{b}_{n}\overline{z}^{n}\right|\\&\geq&
r^{2}_{3}|a_{1}z+\overline{b}_{1}\overline{z}|-r_{3}^{2}\left|\sum_{n=2}^\infty
a_{n}z^{n}+\overline{b}_{n}\overline{z}^{n}\right|\\
&\geq&r^{3}_{3}||a_{1}|-|b_{1}||-r_{3}^{2}\sum_{n=2}^\infty
(|a_{n}|+|b_{n}|)r_{3}^{n}\\
&\geq&r^{3}_{3}-r_{3}^{2}\sum_{n=2}^\infty K(M) r_{3}^{n}\\
&=&r^{3}_{3}- K(M) \cdot\frac{r^{4}_{3}}{1-r_{3}}=\sigma_3.
\end{eqnarray*}

When $M=1$, from Remark 2.5, we get that $g(z)=\alpha z$ or
$g(z)=\alpha\overline{z}$ with $|\alpha|=1$, thus for $|z|=r_{3}=1$,
\begin{eqnarray*}
|F(z)|=|z|^{2}|\alpha ||z|=1.
\end{eqnarray*}
Hence, $F(U_{r_{3}})$ contains the disk $U_{\sigma_{3}}$, where
$r_{3}$ is defined by (\ref{l213}) and $\sigma_{3}$ is defined by
(\ref{l214}).

Finally, it is evident that $r_3=\sigma_3=1$ for $M=1$ is the best
possible. This completes the proof of Theorem 2.10.\hfill $\square$

With the aid of Corollary 2.4, applying the similar method as in our
proof of Theorem 2.10, we can get the following theorem.

{\bf Theorem 2.10$'$} \, Let $g(z)$ be harmonic in the unit disk
$U$, with $g(0)=\lambda_{g}(0)-1=0$ and $|g_{1}(z)|+|g_{2}(z)|\leq
M$ for $z\in U$, as in (\ref{l12}). Then
 $F(z)=|z|^{2}g(z)$ is univalent in the disk $U_{r_{3}'}$, and
$F(U_{r_{3}'})$ contains a schlicht disk $U_{\sigma_{3}'}$, where
\begin{equation}
r_{3}'=\frac{1}{1+2(M-\frac{1}{M})+\sqrt{M-\frac{1}{M}+4(M-\frac{1}{M})^2}},%\label{l213}
\end{equation}
and
\begin{eqnarray}
\sigma_{3}'= \left\{ \begin{array}{lll}
r'^{3}_{3}-(M-\frac{1}{M})\cdot\frac{r'^{4}_{3}}{1-r_{3}'}, &M>1,\\
1 ,\qquad      & M=1.
\end{array} \right.%\label{l214}
\end{eqnarray}
Above result is sharp for $M=1$.

{\bf Remark 2.11} \, Notice that
$M-\frac{1}{M}<K(M)\leq\sqrt{2M^2-2}$ for $M>1$, and
$K(M)=\frac{4M}{\pi}<\sqrt{2M^2-2}$ for
$M>M_0'=\frac{\pi}{\sqrt{\pi^2-8}}\approx 2.2976$, it is easy to
verify that
$$
r_3'>r_3\geq\rho_5,\quad\sigma_{3}'>\sigma_{3}\geq R_5,
$$
for $M>1$.

With the aid of  Theorem C and Lemma 2.1, applying the similar
method as in our proof of Theorem 2.10, we can improve Theorem G as
follows.

{\bf Corollary 2.12} \, Let $g(z)$ be harmonic in the unit disk $U$,
with $g(0)=J_{g}(0)-1=0$ and $|g(z)|\leq M$ for $z\in U$, as in
(\ref{l12}). Then $F(z)=|z|^{2}g(z)$ is univalent in the disk
$U_{r_{4}}$, and $F(U_{r_{4}})$ contains a schlicht disk
$U_{\sigma_{4}}$, where
\begin{equation*}
r_4=\frac{\lambda_0(M)}{\lambda_0(M)+2K(M)+\sqrt{\lambda_0(M)
K(M)+4K(M)^2}},
\end{equation*}
and
\begin{equation*}
\sigma_{4}=\left\{
\begin{array}{lll}
\lambda_0(M) \, r_4^3-K(M)\cdot\frac{r_4^4}{1-r_4}
\, , &\mbox{ if }M>1,\\
1,&\mbox{ if }M=1,
\end{array}
\right.
\end{equation*}
where $\lambda _{0}(M)$ is defined by (\ref{l16}), and $K(M)$ is
defined by (\ref{l220}). Above result is sharp for $M=1$.

With the aid of Lemma 2.3, if we apply the same method as in our
proof of Theorem 2.10, we obtain the fourth version of Landau's
theorem for biharmonic mappings as follows.

{\bf Corollary 2.12$'$} \, Let $g(z)$ be harmonic in the unit disk
$U$, with $g(0)=J_{g}(0)-1=0$ and $|g_{1}(z)|+|g_{2}(z)|\leq M$ for
$z\in U$, as in (\ref{l12}). Then
 $F(z)=|z|^{2}g(z)$ is univalent in the disk $U_{r_{4}'}$, and
$F(U_{r_{4}'})$ contains a schlicht disk $U_{\sigma_{4}'}$, where
\begin{equation*}
r_4'=\frac{\lambda_0(M)}{\lambda_0(M)+2(M-\frac{\lambda_{0}^{2}(M)}{M})+\sqrt{\lambda_0(M)
(M-\frac{\lambda_{0}^{2}(M)}{M})+4(M-\frac{\lambda_{0}^{2}(M)}{M})^2}},
\end{equation*}
and
\begin{equation*}
\sigma_{4}'=\left\{
\begin{array}{lll}
\lambda_0(M) \,
r_4'^3-(M-\frac{\lambda_{0}^{2}(M)}{M})\cdot\frac{r_4'^4}{1-r_4'}
\, , &\mbox{ if }M>1,\\
1,&\mbox{ if }M=1,
\end{array}
\right.
\end{equation*}
where $\lambda _{0}(M)$ is defined by (\ref{l16}). Above result is
sharp for $M=1$.

{\bf Remark 2.13} \, Notice that
$M-\frac{\lambda_{0}^{2}(M)}{M}<K(M)\leq\sqrt{2M^2-2}$ for $M>1$,
and $\lambda_1(M)=\frac{4M}{\pi}<\sqrt{2M^2-2}$ for
$M>M_0'=\frac{\pi}{\sqrt{\pi^2-8}}\approx 2.2976$, it is easy to
verify that
$$
r_4'>r_4>\rho_6>\rho_2,\quad\sigma_{4}'>\sigma_{4}> R_6>R_2,
$$
for $M>M_0'=\frac{\pi}{\sqrt{\pi^2-8}}\approx 2.2976$.

\end{document}